\documentclass[12pt,a4paper]{amsart}
\usepackage[latin1]{inputenc}
\usepackage[T1]{fontenc}
\usepackage[english]{babel}

\usepackage{amsmath,amssymb,amsthm}

\theoremstyle{plain}
\newtheorem{theorem}{Theorem}[section]
\newtheorem{corollary}[theorem]{Corollary}
\newtheorem{prop}[theorem]{Proposition}
\newtheorem{lemma}[theorem]{Lemma}
\theoremstyle{definition}
\newtheorem{remark}[theorem]{Remark}
\newtheorem{example}[theorem]{Example}

\newcommand{\C}{\mathbb{C}}
\newcommand{\D}{\mathbb{D}}
\newcommand{\R}{\mathbb{R}}
\newcommand{\K}{\mathbb{K}}
\newcommand{\N}{\mathbb{N}}
\newcommand{\T}{\mathbb{T}}
\newcommand{\eps}{\varepsilon}

\newcommand{\Id}{\mathrm{Id}}
\newcommand{\en}{\longrightarrow}

\renewcommand{\leq}{\leqslant}
\renewcommand{\geq}{\geqslant}

 \DeclareMathOperator{\re}{Re}
 \DeclareMathOperator{\im}{Im}

\newcounter{equi1}
\newenvironment{equi}
{\begin{list}
        {(\roman{equi1})}
        {\setlength{\itemsep}{0ex plus 0.2ex minus 0ex}
         \setlength{\topsep}{0ex}
         \setlength{\parsep}{0ex}
         \setlength{\labelwidth}{7ex}
         \usecounter{equi1}}}
{\end{list}}

\begin{document}
\title{Norm equalities for operators}
 \subjclass[2000]{Primary: 46B20.}
 \keywords{Daugavet equation; operator norm}
 \date{April 5th, 2006}
 \thanks{The work of the first-named author
was supported by a fellowship from the
\textit{Alexander-von-Humboldt Stiftung}. The work of the
second-named author was partially supported by Spanish MCYT project
no.\ BFM2003-01681 and Junta de Andaluc\'{\i}a grant FQM-185. The work of
the third-named author was supported by a FPI grant of the Spanish
MEC.}

\maketitle

\vspace{0.3cm}

\markboth{\emph{Norm equalities for operators}}{\emph{V.~Kadets,
M.~Mart\'{\i}n, and J.~Mer\'{\i}}}

\centerline{\textsc{\large Vladimir Kadets,}}

\begin{center}\small Faculty of Mechanics and Mathematics \\ Kharkov National
University\\ pl.\ Svobody 4, 61077 Kharkov, Ukraine \\
\emph{E-mail:} \texttt{vova1kadets@yahoo.com}
\end{center}

\vspace{0.3cm}

\centerline{\textsc{\large Miguel Mart\'{\i}n, \qquad and \qquad Javier
Mer\'{\i}}}

\begin{center}\small Departamento de An\'{a}lisis Matem\'{a}tico \\ Facultad de
Ciencias \\ Universidad de Granada \\ 18071 Granada, Spain \\
\emph{E-mail:} \texttt{mmartins@ugr.es, jmeri@ugr.es}
\end{center}

  \thispagestyle{empty}

\begin{abstract}
A Banach space $X$ has the Daugavet property if the Daugavet
equation $\|\Id + T\|= 1 + \|T\|$ holds for every rank-one operator
$T:X \longrightarrow X$. We show that the most natural attempts to
introduce new properties by considering other norm equalities for
operators (like $\|g(T)\|=f(\|T\|)$ for some functions $f$ and $g$)
lead in fact to the Daugavet property of the space. On the other
hand there are equations (for example $\|\Id + T\|= \|\Id - T\|$)
that lead to new, strictly weaker properties of Banach spaces.
\end{abstract}

\section{Introduction}
The purpose of this paper is to study equalities involving the norm
of operators on Banach spaces, and to discuss the possibility of
defining isometric properties for Banach spaces by requiring that
all operators of a suitable class satisfy such a norm equality.

The interest in this topic goes back to 1963, when the Russian
mathematician I.~Daugavet \cite{Dau} showed that each compact
operator $T$ on $C[0,1]$ satisfies the norm equality
\begin{equation*}\label{DE}\tag{DE}
\|\Id + T\|= 1 + \|T\|.
\end{equation*}
The above equation is nowadays referred to as \emph{Daugavet
equation}. Few years later, this result was extended to various
classes of operators on some Banach spaces, including weakly compact
operators on $C(K)$ for perfect $K$ and on $L_1(\mu)$ for atomless
$\mu$ (see \cite{Wer0} for an elementary approach). These results
were left without much attention until the beginning of the
eighties, when a new wave of interest in this topic surfaced, and
the Daugavet equation was studied by many authors in various
contexts. We refer the reader to the books \cite{AbrAli-1,AbrAli-2}
for a brief study of this equation from different points of view.

In the late nineties, new ideas were infused into this field and,
instead of looking for new spaces and new classes of operators on
them for which \eqref{DE} is valid, the geometry of Banach spaces
having the so-called Daugavet property was studied. Following
\cite{KSSW0,KSSW}, we say that a Banach space $X$ has the
\emph{Daugavet property} if every rank-one operator $T\in L(X)$
satisfies \eqref{DE} (we write $L(X)$ for the Banach algebra of all
bounded linear operators on $X$). In such a case, every operator on
$X$ not fixing a copy of $\ell_1$ also satisfies \eqref{DE}
\cite{Shv}; in particular, this happens to every compact or weakly
compact operator on $X$ \cite{KSSW}. There are several
characterizations of the Daugavet property which does not involve
operators (see \cite{KSSW,WerSur}). For instance, a Banach space $X$
has the Daugavet property if and only if for every $x\in S_X$ and
every $\eps>0$ the closed convex hull of the set
$$
B_X\setminus\big(x+(2-\eps)\,B_X\big)
$$
coincides with the whole $B_X$. Here and subsequently, $B_X$ and
$S_X$ stand, respectively, for the closed unit ball and the unit
sphere of a Banach space $X$. Let us observe that the above
characterization shows that the Daugavet property is somehow
extremely opposite to the Radon-Nikod\'{y}m property.

Although the Daugavet property is clearly of isometric nature, it
induces various isomorphic restrictions. For instance, a Banach
space with the Daugavet property does not have the Radon-Nikod\'{y}m
property \cite{Woj} (actually, every slice of the unit ball has
diameter $2$ \cite{KSSW}), it contains $\ell_1$ \cite{KSSW}, it does
not have unconditional basis \cite{Kadets} and, moreover, it does
not isomorphically embed into an unconditional sum of Banach spaces
without a copy of $\ell_1$ \cite{Shv}. It is worthwhile to remark
that the latter result continues a line of generalization
(\cite{K97}, \cite{KS99}, \cite{KSSW}) of the known theorem of
A.~Pe\l czy\'nski \cite{Pelcz} from 1961 saying that neither
$C[0,1]$ nor $L_1[0,1]$ embeds into a space with unconditional basis

The state-of-the-art on the Daugavet property can be found in
\cite{WerSur}; for very recent results we refer the reader to
\cite{BM,BKSW,I-K-W} and references therein.

In view of the deep consequences that the Daugavet property has on
the geometry of a Banach space, one may wonder whether it is
possible to define other interesting properties by requiring all
rank-one operators on a Banach space to satisfy a suitable norm
equality. This is the aim of the present paper.

Let us give some remarks on the question which will also serve to
present the outline of our further discussion. First, the Daugavet
property clearly implies that the norm of $\Id + T$ only depends on
the norm of $T$. Then, a possible generalization of the Daugavet
property is to require that every rank-one operator $T$ on a Banach
space $X$ satisfies a norm equality of the form
$$
\|\Id + T\|=f(\|T\|)
$$
for a fixed function $f:\R^+_0\longrightarrow \R$. It is easy to
show (see Proposition~\ref{prop1}) that the only property which can
be defined in this way is the Daugavet property. Therefore, we
should look for equations in which $\Id + T$ is replaced for other
function of $T$, i.e.\ we fix functions $g$ and $f$ and we require
that every rank-one operator $T$ on a Banach space $X$ satisfies the
norm equality
$$
\|g(T)\|=f(\|T\|).
$$
We need $g$ to carry operators to operators and to apply to
arbitrary rank-one operators, so it is natural to impose $g$ to be a
power series with infinite radius of convergence, i.e.\ an analytic
function defined on the whole $\K$ ($=\R$ or $\C$). Again, the only
non trivial possibility is the Daugavet property, as we will show in
section~\ref{sec:g(T)=f(T)}. Section~\ref{sec:Id+g(T)=f(g(T))} is
devoted to the last kind of equations we would like to study.
Concretely, we consider an analytic function $g$, a continuous
function $f$, and a Banach space $X$, and we require each rank-one
operator $T\in L(X)$ to satisfy the norm equality
\begin{equation}\label{eq:introduction-id+g=fg}
\|\Id + g(T)\|=f(\|g(T)\|).
\end{equation}
If $X$ is a Banach space with the Daugavet property and $g$ is an
analytic function, then it is easy to see that the norm equality
$$
\|\Id + g(T)\|=|1+g(0)|-|g(0)| + \|g(T)\|
$$
holds for every rank-one $T\in L(X)$. Therefore, contrary to the
previous cases, our aim here is not to show that only few functions
$g$ are possible in \eqref{eq:introduction-id+g=fg}, but to prove
that many functions $g$ produce the same property. Unfortunately, we
have to separate the complex case and the real case, and only in the
first one we are able to give fully satisfactory results. More
concretely, we consider a complex Banach space $X$, an entire
function $g$ and a continuous function $f$, such that
\eqref{eq:introduction-id+g=fg} holds for every rank-one operator
$T\in L(X)$. If $\re g(0)\neq -1/2$, then $X$ has the Daugavet
property. Surprisingly, the result is not true when $\re g(0)=
-1/2$, where another family of properties strictly weaker than the
Daugavet property appears: there exists a modulus one complex number
$\omega$ such that the norm equality
\begin{equation}\label{eq:introduction-omega0}
\|\Id + \omega\,T\|=\|\Id+ T\|
\end{equation}
holds for every rank-one $T\in L(X)$. In the real case, the
discussion above depends upon the surjectivity of $g$, and there are
many open questions when $g$ is not onto. Finally, we give in
section~\ref{sec:newproperties} some remarks about the properties
defined by norm equalities of the form given in
\eqref{eq:introduction-omega0}.

We finish the introduction by commenting that, although we have been
unable to find any result as the above ones in the literature, there
are several papers in which the authors work on inequalities which
remind the Daugavet equation. For instance, it is proved in
\cite{Ben-Lin} (see also \cite[\S9]{Pli-Pop}) that for every
$1<p<\infty$, $p\neq 2$, there exists a function
$\varphi_p:(0,\infty)\en (0,\infty)$ such that the inequality
\begin{equation}\label{eq:pseudoDE}
\|\Id + T \| \geq 1 + \varphi_p(\|T\|)
\end{equation}
holds for every nonzero compact operator $T$ on $L_p[0,1]$ (see
\cite{Boy-Kad} for the $\varphi_p$ estimates). This result has been
recently carried to non-commutative $\mathcal{L}_p$-spaces and to
some spaces of operators \cite{Oik-cstar,Oik04}. Finally,
inequalities as \eqref{eq:pseudoDE} in which $\varphi_p$ is linear
are studied in \cite{Mar-daugavetian}.

\section{Preliminaries}\label{sec:preliminares}

Let us start by fixing some notation. We use the symbols $\T$ and
$\D$ to denote, respectively, the unit sphere and the closed unit
ball of the base field $\K$, and we write $\re(\cdot)$ to denote the
real part function, which is nothing than the identity when $\K=\R$.
Given a Banach space $X$, the dual space of $X$ is denoted by
$X^{*}$ and, when $X$ is a complex space, $X_\R$ is the underlying
real space. Given $x\in X$ and $x^*\in X^*$, we write $x^*\otimes x$
to denote the bounded linear operator
$$
y\longmapsto x^*(y)\,x \qquad (y\in X),
$$
whose norm is equal to $\|x\|\,\|x^*\|$.

It is straightforward to check that if an operator $T$ on a Banach
space $X$ satisfies \eqref{DE}, then all the operators of the form
$\alpha\,T$ for $\alpha>0$ also satisfy \eqref{DE} (see
\cite[Lemma~2.1]{AAB} for instance). We will use this fact along the
paper without explicit mention.

Next, let us observe that the definition of the Daugavet property
for complex spaces may be a bit confusing, since the meaning of a
rank-one operator can be understood in two different forms. Indeed,
an operator $T$ on a Banach space $X$ is a rank-one operator if
$\dim T(X)\leq 1$. Then, when $X$ is a complex space, two different
classes of rank-one operators may be considered: the complex-linear
operators whose ranks have complex-dimension~$1$ (i.e.\ operators of
the form $x^*\otimes x$ for $x^*\in X^*$ and $x\in X$) and, on the
other hand, the real-linear operators whose ranks have
real-dimension~$1$ (i.e.\ operators of the form $\re x^*\otimes x$
for $\re x^*\in \left(X_\R\right)^*$ and $x\in X_\R$). As a matter
of facts, the two possible definitions of the Daugavet property that
one may state are equivalent.

\begin{remark}
Let $X$ be a complex Banach space. Then, the following are
equivalent:
\begin{equi}
\item Every (real) rank-one operator on $X_\R$ satisfies
\eqref{DE}.
\item Every (complex) rank-one operator on $X$ satisfies \eqref{DE}.
\end{equi}
\end{remark}

\begin{proof}
$(i)\Rightarrow (ii)$ follows immediately from
\cite[Theorem~2.3]{KSSW}, but a direct proof is very easy to state.
Given a complex rank-one operator $T=x_0^*\otimes x_0$ with
$\|x_0^*\|=1$ and $\|x_0\|=1$, we have $\|\Id+\re x_0^*\otimes
x_0\|=2$. Therefore, given $\eps>0$, there exists $x\in S_X$ such
that $\|x+\re x_0^*(x)x_0\|\geq 2-\eps$. It follows that
$$
|\re x_0^*(x)|\geq 1-\eps, \qquad \text{and so} \qquad |\im
x_0^*(x)|\leq \eps.
$$
Now, it is clear that
\begin{align*}
\|\Id+T\|&\geq \|x + T x\| =\|x+x_0^*(x)x_0\| \\ &\geq \|x+\re
x_0^*(x)x_0\| -\|\im x_0^*(x)x_0\|\geq 2-2\eps.
\end{align*}
$(ii)\Rightarrow (i)$. Let $T=\re x_0^* \otimes x_0$ be a real
rank-one operator with
$$
\|\re x_0^*\|=\|x_0^*\|=1 \qquad \text{and} \qquad \|x_0\|=1.
$$
By $(ii)$, we have that $\|\Id+x_0^*\otimes x_0\|=2$ so, given
$\eps>0$, there exists $x\in S_X$ satisfying
$$
\|x+x_0^*(x)x_0\|\geq2-\eps.
$$
If we take $\omega\in\T$ such that $\omega x_0^*(x)=|x_0^*(x)|$,
then
$$
x_0^*(\omega x)=\omega x_0^*(x)=\re x_0^*(\omega x).
$$
Therefore,
\begin{align*}
\|\Id+T\| &\geq \|\omega x+\re x_0^*(\omega x)x_0\| \\ &= \|\omega
x+ \omega x_0^*(x)x_0\|= \|x+x_0^*(x)x_0\|\geq 2-\eps.\qedhere
\end{align*}
\end{proof}

From now on, by a rank-one operator on a Banach space over $\K$, we
will mean a bounded $\K$-linear operator whose image has
$\K$-dimension less or equal than $1$.

As we commented in the introduction, the aim of this paper is to
discuss whether there are other isometric properties apart from the
Daugavet property which can be defined by requiring all rank-one
operators on a Banach space to satisfy a norm equality. A first
observation is that the Daugavet property implies that for every
rank-one operator $T$, the norm of $\Id + T$ only depends on the
norm of $T$. It is easy to check that the above fact only may happen
if the Banach space involved has the Daugavet property. We state and
prove a slightly more general version of this result which we will
use later on.

\begin{prop}\label{prop1}
Let $f:\R^+_0\longrightarrow\R^+_0$ be an arbitrary function.
Suppose that there exist $a, b\in \K$ and a non-null Banach space
$X$ over $\K$ such that the norm equality
$$
\|a\Id+b\,T\|=f(\|T\|)
$$
holds for every rank-one operator $T\in L(X)$. Then,
$f(t)=|a|+|b|\,t$ for every $t\in\R^+_0$. In particular, if $a\neq
0$ and $b\neq 0$, then $X$ has the Daugavet property.
\end{prop}

\begin{proof}
If $a\,b=0$ we are trivially done, so we may assume that $a\neq0$,
$b\neq0$ and we write
$\omega_0=\frac{\overline{b}}{|b|}\frac{a}{|a|}\in \T$. Now, we fix
$x_0\in S_X$, $x_0^*\in S_{X^*}$ such that $x_0^*(x_0)=\omega_0$
and, for each $t\in\R^+_0$, we consider the rank-one operator
$T_t=t\,x_0^* \otimes x_0\in L(X)$. Observe that $\|T_t\|=t$, so we
have
$$
f(t)=\|a \Id+ b\, T_t\| \qquad (t\in \R^+_0).
$$
Then, it follows that
\begin{align*}
|a|+|b|\,t \geq f(t) &= \|a \Id + b\, T_t\|\geq\bigl\|[a\Id+b\,
T_t](x_t)\bigr\| \\ & = \|a\,x_t+b\,\omega_0t\,x_t\|=|a +
b\,\omega_0 t|\,\|x_t\|\\ & =\left|a + b
\frac{\overline{b}}{|b|}\frac{a}{|a|}\, t \right| = |a|+|b|\,t.
\end{align*}
Finally, if the norm equality
$$
\|a \Id + b\, T\|= |a| + |b|\, \|T\|
$$
holds for every rank-one operator on $X$, then $X$ has the Daugavet
property. Indeed, we fix a rank-one operator $T\in L(X)$ and apply
the above equality to $S=\frac{a}{b}\,T$ to get
\begin{equation*}
|a|\bigl(1 + \|T\|\bigr)=|a|+ |b|\,\|S\| =\|a\Id + b\, S\|= |a|\,
\|\Id +T\|.\qedhere
\end{equation*}
\end{proof}

With the above property in mind, we have to look for Daugavet-type
norm equalities in which $\Id + T$ is replaced by other function of
$T$. If we want such a function to carry operators to operators and
to be applied to arbitrary rank-one operators on arbitrary Banach
spaces, it is natural to consider $g$ to be a power series with
infinite radius of convergence. In other words, we consider an
analytic function $g:\K\en \K$ and, for each operator $T\in L(X)$,
we define
$$
g(T)=\sum_{k=0}^\infty a_k\, T^k,
$$
where $g(\zeta)=\sum_{k=0}^\infty a_k\, \zeta^k$ is the power series
expansion of $g$. The following easy result shows how to calculate
$g(T)$ when $T$ is a rank-one operator.

\begin{lemma}\label{lemma}
Let $g:\K \en \K$ be an analytic function with power series
expansion
$$
g(\zeta)=\sum_{k=0}^\infty a_k\,\zeta^k \qquad (\zeta\in \K),
$$
and let $X$ be a Banach space over $\K$. For $x^*\in X^*$ and $x\in
X$, we write $T=x^*\otimes x$ and $\alpha=x^*(x)$. Then, for each
$\lambda\in\K$,
$$
g(\lambda T)=\begin{cases} a_0 \Id + a_1\lambda T & \text{if $\alpha=0$} \\
a_0 \Id + \frac{\widetilde{g}(\alpha\lambda)}{\alpha}T & \text{if
$\alpha\neq 0$,}
\end{cases}
$$
where
$$
\widetilde{g}(\zeta)=g(\zeta)-a_0 \qquad \bigl(\zeta\in \K\bigr).
$$
\end{lemma}

\begin{proof}
Given $\lambda\in\K$, it is immediate to check that
$$
(\lambda T)^k=\alpha^{k-1}\lambda^k\,T \qquad (k\in\N).
$$
Now, if $\alpha=0$, then $T^2=0$ and the result is clear. Otherwise,
we have
\begin{align*}
g(\lambda T) & =a_0 \Id + \sum_{k=1}^\infty
a_k\,\alpha^{k-1}\,\lambda^k\,T \\ & = a_0 \Id +
\left(\frac{1}{\alpha}\sum_{k=1}^\infty
a_k\,\alpha^k\lambda^k\,\right)T = a_0 \Id +
\frac{\widetilde{g}(\alpha\lambda)}{\alpha}T.\qedhere
\end{align*}
\end{proof}

\section{Norm equalities of the form
$\|g(T)\|=f(\|T\|)$} \label{sec:g(T)=f(T)}

We would like to study now norm equalities for operators of the form
\begin{equation}\label{eq:normg(T)=f(normT)}
\|g(T)\|=f(\|T\|),
\end{equation}
where $f:\R^+_0\en \R^+_0$ is an arbitrary function and $g:\K\en \K$
is an analytic function.

Our target is to show that the Daugavet property is the only
non-trivial property that it is possible to define by requiring all
rank-one operators on a Banach space of dimension greater than one
to satisfy a norm equality of the form \eqref{eq:normg(T)=f(normT)}.
We start by proving that $g$ has to be a polynomial of degree less
or equal than $1$, and then we will deduce the result from
Proposition~\ref{prop1}.

\begin{theorem}\label{theorem}
Let $g:\K\longrightarrow \K$ be an analytic function and $f: \R_0^+
\longrightarrow \R_0^+$ an arbitrary function. Suppose that there is
a Banach space $X$ over $\K$ with $\dim(X)\geq 2$ such that the norm
equality
\begin{equation*}
\|g(T)\|=f(\|T\|)
\end{equation*}
holds for every rank-one operator $T$ on $X$. Then, there are $a,\,
b\in \K$ such that
$$
g(\zeta)=a+b\zeta \qquad \bigl(\zeta \in \K\bigr).
$$
\end{theorem}

\begin{proof} Let $g(\zeta)= \sum_{k=0}^\infty a_k\, \zeta^k$ be the
power series expansion of $g$ and let \mbox{$\widetilde{g}=g-a_0$}.
Given $\alpha\in \D$, we take $x_\alpha^*\in S_{X^*}$ and
$x_\alpha\in S_X$ such that $x_\alpha^*(x_\alpha)=\alpha$ (we can do
it since $\dim(X)\geq2$), and we write $T_\alpha=x_\alpha^*\otimes
x_\alpha$, which satisfies $\|T_\alpha\|=1$. Using
Lemma~\ref{lemma}, for each $\lambda\in \K$ we obtain that
\begin{align*}
g(\lambda T_0) &=a_0\Id+a_1\lambda\, T_0 \intertext{and} g(\lambda
T_\alpha) &=
a_0\Id+\frac{1}{\alpha}\widetilde{g}(\lambda\alpha)\,T_\alpha \qquad
(\alpha\neq 0).
\end{align*}
Now, fixed $\lambda\in \K$, we have
\begin{align*}
f(|\lambda|)&=\|g(\lambda T_0)\|=\|a_0\Id +a_1\lambda T_0\|,
\intertext{and} f(|\lambda|)&=\|g(\lambda T_\alpha)\|=
\left\|a_0\Id+\frac{1}{\alpha}\widetilde{g}(\lambda\alpha)T_\alpha\right\|.
\end{align*}
Therefore, we have the equality
\begin{equation}\label{eq:opzero=opalpha}
\left\|a_0\Id+\frac{1}{\alpha}\widetilde{g}(\lambda\alpha)T_\alpha\right\|=\|a_0\Id
+a_1\lambda T_0\|\qquad (\lambda\in\K,\ 0<|\alpha|\leq1).
\end{equation}
In the complex case is enough to consider the above equality for
$\alpha=1$ and to use the triangle inequality to get that
\begin{equation}\label{eq:th1-complex}
|\widetilde{g}(\lambda)|\leq 2|a_0| + |a_1|\,|\lambda| \qquad
(\lambda\in \C).
\end{equation}
From this, it follows by just using Cauchy's estimates, that
$\widetilde{g}$ is a polynomial of degree one (see \cite[Exercise 1,
p.~80]{Conway} or \cite[Theorem~3.4.4]{Gre-Kra}, for instance), and
we are done.

In the real case, it is not possible to deduce from inequality
\eqref{eq:th1-complex} that $\widetilde{g}$ is a polynomial, so we
have to return to \eqref{eq:opzero=opalpha}. From this equality, we
can deduce by just applying the triangle inequality that
\begin{equation*}
\left|\frac{\widetilde{g}(\lambda\alpha)}{\alpha}\right|-|a_0|\leq|a_0|+|a_1|\,|\lambda|
\qquad \text{and} \qquad
|a_1|\,|\lambda|-|a_0|\leq\left|\frac{\widetilde{g}(\lambda\alpha)}{\alpha}\right|+|a_0|
\end{equation*}
for every $\lambda\in\R$ and every $\alpha\in [-1,1]\setminus\{0\}$.
It follows that
\begin{equation}\label{th:eq-desigualdad-importante}
\left|\left|\frac{\widetilde{g}(\lambda\alpha)}{\alpha}\right|-|a_1|\,|\lambda|\right|\leq
2|a_0|\qquad \bigl(\lambda\in\R,\ \alpha\in
[-1,1]\setminus\{0\}\bigr).
\end{equation}
Next, for $t\in ]1,+\infty[$ and $k\in\N$, we use
\eqref{th:eq-desigualdad-importante} with $\lambda=t^k$ and
$\alpha=\frac{1}{t^{k-1}}$ to obtain that
$$
\big|\,|\widetilde{g}(t)|-|a_1|\,t\big|\leq\frac{2|a_0|}{t^{k-1}}
$$
so, letting $k\longrightarrow\infty$, we get that
$$
|\widetilde{g}(t)|=|a_1|\,t \qquad \bigl(t\in ]1,+\infty[\bigr).
$$
Finally, an obvious continuity argument allows us to deduce from the
above equality that $\widetilde{g}$ coincides with a degree one
polynomial in the interval $]1,+\infty[$, thus the same is true in
the whole $\R$ by analyticity.
\end{proof}

We summarize the information given in Proposition~\ref{prop1} and
Theorem~\ref{theorem}.

\begin{corollary}
Let $f: \R_0^+ \longrightarrow \R_0^+$ be an arbitrary function and
$g:\K\longrightarrow \K$ an analytic function. Suppose that there is
a Banach space $X$ over $\K$ with $\dim(X)\geq 2$ such that the norm
equality
\begin{equation*}
\|g(T)\|=f(\|T\|)
\end{equation*}
holds for every rank-one operator $T$ on $X$. Then, only three
possibilities may happen:
\begin{enumerate}
\item[(a)] $g$ is a constant function (trivial case).
\item[(b)] There is a non-null $b\in \K$ such that
$g(\zeta)=b\,\zeta$ for every $\zeta\in \K$ (trivial case).
\item[(c)] There are non-null $a,b\in \K$ such that $g(\zeta)=a +
b\, \zeta$ for every $\zeta\in \K$, and $X$ has the Daugavet
property.
\end{enumerate}
\end{corollary}

\begin{remark}
{\slshape The above result is not true when the dimension of the
space is one.\ } Indeed, for every $a\in L(\K)\equiv \K$, one has
$|a^2|=|a|^2$.
\end{remark}

\section{Norm equalities of the form $\|\Id + g(T)\|=f(\|g(T)\|)$}
\label{sec:Id+g(T)=f(g(T))}

Let $X$ be a Banach space over $\K$. Our next aim is to study norm
equalities of the form
\begin{equation}\label{eq:DEwith-g-and-f(g)}
\|\Id + g(T)\| = f(\|g(T)\|)
\end{equation}
where $g:\K \en \K$ is analytic and $f:\R^+_0 \en \R^+_0$ is
continuous.

When $X$ has the Daugavet property, it is clear that
Eq.~\eqref{eq:DEwith-g-and-f(g)} holds for every rank-one operator
if we take $g(\zeta)=\zeta$ and $f(t)=1 + t$. But, actually, every
analytic function $g$ works with a suitable $f$.

\begin{remark}\label{remark:4.1}
{\slshape If $X$ is a real or complex Banach space with the Daugavet
property and $g:\K \en \K$ is an analytic function, the norm
equality
$$
\|\Id + g(T)\|=|1+g(0)|-|g(0)| + \|g(T)\|
$$
holds for every weakly compact operator $T\in L(X)$.\ } Indeed,
write $\widetilde{g}=g-g(0)$ and observe that $\widetilde{g}(T)$ is
weakly compact whenever $T$ is. Since $X$ satisfies the Daugavet
property, we have
\begin{align*}
\|\Id + g(T)\| &=\|\Id + g(0)\,\Id + \widetilde{g}(T)\|=|1+g(0)| +
\|\widetilde{g}(T)\| \\ & = \bigl(|1+g(0)| -|g(0)|\bigr) +
\bigl(|g(0)| + \|\widetilde{g}(T)\|\bigr) \\ & =|1+g(0)|-|g(0)| +
\|g(0)\,\Id + \widetilde{g}(T)\| \\ &= |1+g(0)|-|g(0)| + \|g(T)\|.
\end{align*}
\end{remark}

With the above result in mind, it is clear that the aim of this
section cannot be to show that only few $g$'s are possible in
\eqref{eq:DEwith-g-and-f(g)}, but it is to show that many $g$'s
produce only few properties. At first, before formulating our
results, let us discuss the case when the Banach space we consider
is one-dimensional.

\begin{remark}$ $\label{remark-dim-one}
\begin{enumerate}
\item[(a)] \textsc{Complex case:\ } {\slshape It is not possible
to find a non-constant entire function $g$ and an arbitrary function
$f:\R^+_0\en \R$ such that the equality
$$
|1+g(\zeta)|=f(|g(\zeta)|)
$$
holds for every $\zeta\in \C\equiv L(\C)$.\ } Indeed, we suppose
otherwise that such functions $g$ and $f$ exist, and we use Picard
Theorem to assure the existence of $\lambda>0$ such that
$-\lambda,\, \lambda \in g(\C)$. We get
$$
f(\lambda)=|1+\lambda|\neq |1-\lambda|=f(|-\lambda|)=f(\lambda),
$$
a contradiction.
\item[(b)] \textsc{Real case:\ } {\slshape The equality
$$
|1+t^2|=1 + |t^2|
$$
holds for every $t\in \R\equiv L(\R)$.}
\end{enumerate}
\end{remark}

It follows that real and complex spaces do not behave in the same
way with respect to equalities of the form given in
\eqref{eq:DEwith-g-and-f(g)}. Therefore, from now on we study
separately the complex and the real cases. Let us also remark that
when a Banach space $X$ has dimension greater than one, it is clear
that
$$
\|g(T)\|\geqslant |g(0)|
$$
for every analytic function $g:\K\en \K$ and every rank-one operator
$T\in L(X)$. Therefore, the function $f$ in
\eqref{eq:DEwith-g-and-f(g)} has to be defined only in the interval
$[|g(0)|,+\infty[$.

\vspace{0.5cm}

\noindent \textbf{\large $\bullet$} \textsc{\large Complex case:}

Our key lemma here states that the function $g$ in
\eqref{eq:DEwith-g-and-f(g)} can be replaced by a degree one
polynomial.

\begin{lemma}\label{lemma:4.3}
Let $g:\C\en\C$ be a non-constant entire function, let\linebreak
$f:[|g(0)|,+\infty[\en \R$ be a continuous function and let $X$ be a
Banach space with dimension greater than $1$. Suppose that the norm
equality
$$
\|\Id + g(T)\| = f(\|g(T)\|)
$$
holds for every rank-one operator $T\in L(X)$. Then,
$$
\bigl\|\big(1+g(0)\big)\,\Id + T\bigr\|=|1+g(0)|-|g(0)| +
\|g(0)\,\Id + T\|
$$
for every rank-one operator $T\in L(X)$.
\end{lemma}

\begin{proof}
We claim that the norm equality
\begin{equation}\label{eq:gtopoldegree1}
\bigl\|\big(1+g(0)\big)\,\Id + T\bigr\|=f\big(\|g(0)\,\Id + T\|\big)
\end{equation}
holds for every rank-one operator $T\in L(X)$. Indeed, we write
$\widetilde{g}=g-g(0)$ and we use Picard Theorem to assure that
$\widetilde{g}(\C)$ is equal to $\C$ except for, eventually, one
point $\alpha_0\in \C$. Next, we fix a rank-one operator
$T=x^*\otimes x$ with $x^*\in X^*$ and $x\in X$. If $x^*(x)\neq 0$
and $x^*(x)\neq \alpha_0$, we may find $\zeta\in \C$ such that
$\widetilde{g}(\zeta)=x^*(x)$, and we use Lemma~\ref{lemma} to get
that
$$
g\left(\frac{\zeta}{x^*(x)}\,T\right)= g(0)\,\Id +
\frac{\widetilde{g}(\zeta)}{x^*(x)}\,T=g(0)\,\Id +T.
$$
We deduce that
\begin{align*}
\left\|\bigl(1 + g(0)\bigr)\,\Id + T\right\|& =\left\|\Id +
g\left(\frac{\zeta}{x^*(x)}\,T\right)\right\| \\ & =
f\left(\left\|g\left(\frac{\zeta}{x^*(x)}\,T\right)\right\|\right) =
f\bigl(\|g(0)\,\Id + T\|\bigr).
\end{align*}
The remaining cases in which $x^*(x)=0$ or $x^*(x)=\alpha_0$ follow
from the above equality thanks to the continuity of $f$.

To finish the proof, we have to show that
\begin{equation*}
f(t)=|1+g(0)|-|g(0)| + t \qquad \bigl(t\geq |g(0)|\bigr).
\end{equation*}
Suppose first that $g(0)=-1$. We take $x\in S_X$ and $x^*\in
S_{X^*}$ such that \mbox{$x^*(x)=1$} and, for every $t\geq 1$, we
define the rank-one operator
$$
T_t=(1-t)\,x^*\otimes x.
$$
It is immediate to show that
$$
\|-\Id + T_t\|=t \qquad \text{and} \qquad \|T_t\|=t-1.
$$
Then, it follows from \eqref{eq:gtopoldegree1} that $f(t)=t-1$, and
we are done.

Suppose otherwise that $g(0)\neq -1$. We take $x\in S_X$ and $x^*\in
S_{X^*}$ such that $x^*(x)=1$ and, for every $t\geq |g(0)|$, we
define the rank-one operator
$$
T_t=\frac{1+g(0)}{|1+g(0)|}\,\big(t-|g(0)|\big)\,x^*\otimes x.
$$
It is routine to show, by just evaluating at the point $x$, that
\begin{equation*}
\left\|\bigl(1 + g(0)\bigr)\,\Id + T_t\right\|=|1+g(0)|+t-|g(0)|.
\end{equation*}
Therefore, if follows from \eqref{eq:gtopoldegree1} that
\begin{equation}\label{eq:gtopoldegree1(2)}
f\bigl(\|g(0)\,\Id + T_t\|\bigr)=|1+g(0)|+t-|g(0)|,
\end{equation}
and the proof finishes by just proving that
$$
\|g(0)\,\Id+T_t\|=t.
$$
On the one hand, it is clear that
$$
\|g(0)\,\Id+T_t\|\leq |g(0)| + \|T_t\|=|g(0)| + t - |g(0)|=t.
$$
On the other hand, the converse inequality trivially holds when
$g(0)=0$, so we may suppose \mbox{$g(0)\neq 0$} and we define the
rank-one operator
$$
S_t=\frac{g(0)}{|g(0)|}\,\bigl(\|g(0)\,\Id+T_t\|-|g(0)|\bigr)\,x^*\otimes
x.
$$
It is routine to show, by using that $\|g(0)\,\Id+S_t\|\geq |g(0)|$
and evaluating at the point $x$, that
$$
\|g(0)\,\Id+S_t\|=\|g(0)\,\Id+T_t\|.
$$
From the above equality, \eqref{eq:gtopoldegree1}, and
\eqref{eq:gtopoldegree1(2)}, we deduce that
\begin{multline*}
|1+g(0)|+t-|g(0)| = f\bigl(\|g(0)\,\Id +
S_t\|\bigr)=\left\|\bigl(1+g(0)\bigr)\,\Id+S_t\right\| \\
\leq |1+g(0)| + \|S_t\| =|1+g(0)|+\|g(0)\,\Id + T_t\|-|g(0)|,
\end{multline*}
so $t\leq \|g(0)\,\Id+T_t\|$ and we are done.
\end{proof}

In view of the norm equality appearing in the above lemma, two
different cases arise either $|1+g(0)|\neq |g(0)|$ or $|1+g(0)|=
|g(0)|$; equivalently, $\re g(0)\neq -1/2$ or $\re g(0)=-1/2$. In
the first case, we get the Daugavet property.

\begin{theorem}\label{th:4.4}
Let $X$ be a complex Banach space with $\dim(X)\geq2$. Suppose that
there exist a non-constant entire function $g:\C\longrightarrow\C$
with $\re g(0)\neq-\frac{1}{2}$ and a continuous function
$f:[|g(0)|,+\infty[\longrightarrow\R^+_0$, such that the norm
equality
$$
\|\Id+g(T)\|=f(\|g(T)\|)
$$
holds for every rank-one operator $T\in L(X)$. Then, $X$ has the
Daugavet property.
\end{theorem}

\begin{proof}
Let us first suppose that $\re g(0)>-1/2$ so that
$$
M=|1+g(0)|-|g(0)|>0.
$$
Then, dividing by $M$ the equation given by Lemma~\ref{lemma:4.3},
we get that the norm equality
\begin{equation}\label{eq:reg(0)>-1/2---1}
\left\|\textstyle\frac{1+g(0)}{M}\,\Id+\frac{1}{M}\,T\right\|=1+\left\|\textstyle
\frac{g(0)}{M}\,\Id+ \frac{1}{M}\,T\right\|
\end{equation}
holds for every rank-one operator $T\in L(X)$. Now, we take
$\omega,\,\xi\in\T$ such that
$$
\omega\,(1+g(0))=|1+g(0)|\qquad \text{and} \qquad \xi\, g(0)=|g(0)|,
$$
we fix a rank-one operator $T\in L(X)$, and we observe that
$$
\left\|\textstyle{\frac{1+g(0)}{M}}\,\Id+\frac{1}{M}\,T\right\|=
\left\|\textstyle{\frac{|1+g(0)|}{M}}\,\Id+\frac{\omega}{M}\,T\right\|=
\left\|\Id+\textstyle{\frac{|g(0)|}{M}}\,\Id+\frac{\omega}{M}\,T\right\|
$$
and
$$
\left\|\textstyle\frac{g(0)}{M}\,\Id+ \frac{1}{M}\,T\right\|=
\left\|\textstyle\frac{|g(0)|}{M}\,\Id+\frac{\xi}{M}\,T\right\|.
$$
Therefore, from \eqref{eq:reg(0)>-1/2---1} we get
$$
\left\|\textstyle\Id+\frac{|g(0)|}{M}\,\Id+\frac{\omega}{M}\,T\right\|=
1+\left\|\textstyle\frac{|g(0)|}{M}\,\Id+\frac{\xi}{M}\,T\right\|.
$$
It follows straightforwardly from the arbitrariness of $T$ that the
norm equality
\begin{equation}\label{eq:teor-complejo-igualdad-xi-omega}
\left\|\textstyle\Id+\frac{|g(0)|}{M}\,\Id+T\right\|=
1+\left\|\textstyle\frac{|g(0)|}{M}\,\Id+\frac{\xi}{\omega}\,T\right\|
\end{equation}
holds for every rank-one operator $T\in L(X)$. Now, we claim that
\begin{equation}\label{eq:teor-complejo-----2}
\left\|\textstyle\frac{|g(0)|}{M}\,\Id+T\right\|=
\left\|\textstyle\frac{|g(0)|}{M}\,\Id+\frac{\xi}{\omega}\,T\right\|
\end{equation}
for every rank-one operator $T\in L(X)$. Indeed, by the triangle
inequality, we deduce from
\eqref{eq:teor-complejo-igualdad-xi-omega} that
$$
\left\|\textstyle\frac{|g(0)|}{M}\,\Id+T\right\|\geq
\left\|\textstyle\frac{|g(0)|}{M}\,\Id+\frac{\xi}{\omega}\,T\right\|
$$
for every rank-one operator $T\in L(X)$ and, for every $n\in \N$,
applying the above inequality $n$ times, we get that
$$
\left\|\textstyle\frac{|g(0)|}{M}\,\Id+T\right\|\geq
\left\|\textstyle\frac{|g(0)|}{M}\,\Id+\frac{\xi}{\omega}\,T\right\|\geq
\cdots \geq
\left\|\textstyle\frac{|g(0)|}{M}\,\Id+\left(\frac{\xi}{\omega}\right)^n\,T\right\|;
$$
the claim follows from the easy fact that the sequence
$\left\{\left(\frac{\xi}{\omega}\right)^n\right\}_{n\in\N}$ has a
subsequence which converges to $1$.

Now, given a rank-one operator $T\in L(X)$, it follows from
\eqref{eq:teor-complejo-igualdad-xi-omega} and
\eqref{eq:teor-complejo-----2} that
$$
\left\|\textstyle\Id+\frac{|g(0)|}{M}\,\Id+T\right\|= 1+
\left\|\textstyle\ \frac{|g(0)|}{M}\,\Id+T\right\|
$$
holds and, therefore, for every $n\in\N$ we get
$$
\left\|\textstyle\Id+\frac{1}{n}\left(\textstyle\frac{|g(0)|}{M}\,\Id+T\right)\right\|=
1 + \left\|\textstyle\frac{1}{n}
\left(\textstyle\frac{|g(0)|}{M}\,\Id+T\right)\right\|.
$$
To finish the proof, we fix a rank-one operator $S\in L(X)$ and
$n\in\N$, and we apply the above equality to $T=n\,S$ to get that
$$
\left\|\textstyle\Id+S+\textstyle\frac{|g(0)|}{n\,M}\,\Id\right\|=
\left\|\textstyle\frac{|g(0)|}{n\,M}\,\Id+S\right\|.
$$
We let $n\to \infty$ to deduce that
\begin{equation*}
\|\Id+S\|=1+\|S\|,
\end{equation*}
so $X$ has the Daugavet property.

In case that $\re g(0)<-\frac{1}{2}$, we write
$$
M=-|1+g(0)|+|g(0)|>0
$$
and we deduce from Lemma~\ref{lemma:4.3} that the norm equality
\begin{equation*}
\left\|\textstyle \frac{g(0)}{M}\,\Id+
\frac{1}{M}\,T\right\|=1+\left\|\textstyle\frac{1+g(0)}{M}\,\Id+\frac{1}{M}\,T\right\|
\end{equation*}
holds for every rank-one operator $T\in L(X)$. The rest of the proof
is completely analogous, using the above equality instead of
\eqref{eq:reg(0)>-1/2---1}.
\end{proof}

When $\re g(0)=-\frac{1}{2}$, the above proof does not work.
Actually, contrary to all the previous cases, another family of
properties apart from the Daugavet property appears. More
concretely, let $X$ be a complex Banach space with dimension greater
than two, let $g:\C\en\C$ be a non-constant entire function with
$\re g(0)=-1/2$, and let $f:[|g(0)|,+\infty[\en \R$ be a continuous
function, such that the norm equality
$$
\|\Id + g(T)\|=f(\|g(T)\|)
$$
holds for every rank-one operator $T\in L(X)$. Then, we get from
Lemma~\ref{lemma:4.3} that
$$
\bigl\|\big(1+g(0)\big)\,\Id + T\bigr\|=\|g(0)\,\Id + T\|
$$
for every rank-one operator $T\in L(X)$. Therefore, being
$|1+g(0)|=|g(0)|$, we deduce that there are $\omega_1,\,\omega_2\in
\C$ with $\omega_1\neq\omega_2$ and $|\omega_1|=|\omega_2|$, such
that
$$
\|\Id+\omega_1\,T\|=\|\Id+\omega_2\,T\|
$$
for every rank-one operator $T\in L(X)$ or, equivalently, that there
is $\omega\in\T\setminus\{1\}$ such that
\begin{equation*}
\|\Id + \omega\,T\|=\|\Id+T\|
\end{equation*}
for every rank-one operator $T\in L(X)$. It is routine to check
that, fixed a Banach space $X$, the set of those $\omega\in \T$
which make true the above equality for all rank-one operators on $X$
is a multiplicative closed subgroup of $\T$. Recall that such a
subgroup of $\T$ is either the whole $\T$ or the set of those
$n^{\rm th}$-roots of unity for a positive entire $n\geq 2$. Let us
state the result the we have just proved.

\begin{theorem}\label{th:4.5}
Let $X$ be a complex Banach space with $\dim(X)\geq2$. Suppose that
there exist a non-constant entire function $g:\C\longrightarrow\C$
with $\re g(0)=-\frac{1}{2}$ and a continuous function
$f:[|g(0),+\infty[\longrightarrow\R^+_0$, such that the norm
equality
$$
\|\Id+g(T)\|=f(\|g(T)\|)
$$
holds for every rank-one operator $T\in L(X)$. Then, there is
$\omega\in \T\setminus\{1\}$ such that
$$
\|\Id + \omega T\|=\|\Id + T\|
$$
for every rank-one operator $T\in L(X)$. Moreover, two possibilities
may happen:
\begin{enumerate}
\item[(a)] If $\omega^n\neq 1$ for every $n\in \N$, then
$$
\|\Id + \xi\,T\|=\|\Id + T\|
$$
for every rank-one operator $T\in L(X)$ and every $\xi\in \T$.
\item[(b)] Otherwise, if we take the minimum $n\in\N$ such that
$\omega^n=1$, then
$$
\|\Id + \xi\,T\|=\|\Id + T\|
$$
for every rank-one operator $T\in L(X)$ and every $n^{\rm th}$-root
$\xi$ of unity.
\end{enumerate}
\end{theorem}

The next example shows that all the properties appearing above are
strictly weaker than the Daugavet property.

\begin{example}\label{example-new-property}
{\slshape The real or complex Banach space $X=C[0,1]\oplus_2C[0,1]$
does not have the Daugavet property. However, the norm equality
\begin{equation*}
\|\Id+\omega T\|=\|\Id+T\|
\end{equation*}
holds for every rank-one operator $T\in L(X)$ and every
$\omega\in\T$.}
\end{example}

\begin{proof}
$X$ does not have the Daugavet property by
\cite[Corollary~5.4]{BKSW}. For the second assertion, we fix
$\omega\in\T$ and a rank-one operator $T=x^*\otimes x$ on $X$ (we
take $x\in S_X$ and $x^*\in X^*$), and it is enough to check that
$$
\|\Id + \omega\,T\|^2\geq \|\Id + T\|^2.
$$
To do so, we write $x^*=(x_1^*,x_2^*)$ and $x=(x_1,x_2)$, with
$x_1^*$, $x_2^*\in C[0,1]^*$ and $x_1$, $x_2\in C[0,1]$, and we may
and do assume that
$$
x_1^*=\mu_1+\sum_{j=1}^{n_1}\alpha_j\delta_{r_j}\qquad
x_1^*=\mu_2+\sum_{j=1}^{n_2}\beta_j\delta_{s_j},
$$
where $\alpha_1$,$\dots$, $\alpha_{n_1}$, $\beta_1$,$\dots$,
$\beta_{n_2}\in \C$, $r_1$,$\dots$, $r_{n_1}$, $s_1$,$\dots$,
$s_{n_2}\in [0,1]$, and $\mu_1$, $\mu_2$ are non-atomic measures on
$[0,1]$ (indeed, each rank-one operator can be approximated by
operators satisfying the preceding condition).

Now, we fix $0<\eps<1$ and we consider $y=(y_1,y_2)\in S_X$ such
that
\begin{equation}\label{example:norm(Id+T)}
\|y+Ty\|^2=\|y_1+x^*(y)x_1\|^2+\|y_2+x^*(y)x_2\|^2\geq\|\Id+T\|^2-\eps.
\end{equation}
Since $x_1$, $x_2$, $y_1$, and $y_2$ are continuous functions and
$[0,1]$ is perfect, we can find open intervals $\Delta_1$,
$\Delta_2\subset [0,1]$ so that
\begin{equation}\label{open-interval-property}
|y_i(t)+x^*(y)x_i(t)|\geq\|y_i+x^*(y)x_i\|-\eps \qquad \big(t\in
\Delta_i,\ i=1,2\big)
\end{equation}
and
\begin{equation}\label{open-interval-property2}
\Delta_1\cap\{r_j \ : \ j=1,\dots,n_1\}=\emptyset, \qquad
\Delta_2\cap\{s_j \ : \ j=1,\dots, n_2\}=\emptyset
\end{equation}
Furthermore, using that $\mu_1$, $\mu_2$ are non-atomic and reducing
$\Delta_1$, $\Delta_2$ if necessary, we can assume that they also
satisfy
\begin{equation}\label{open-interval-property2-bis}
|\mu_1(f)|<\eps\,\|f\|, \qquad |\mu_2(g)|<\eps\,\|f\|
\end{equation}
for every $f,g\in C[0,1]$ with $\text{supp}(f)\subset \Delta_1$,
$\text{supp}(g)\subset \Delta_2$.

Now, we fix $t_1\in\Delta_1$ and $t_2\in\Delta_2$. For $i=1,2$, we
take piece-wise linear continuous functions
$\phi_i:[0,1]\longrightarrow [0,1]$ such that
$$
\phi_i(t_i)=1 \qquad \text{and} \qquad
\phi_i\big([0,1]\setminus\Delta_i\big)=\{0\},
$$
and we define
$$
\widetilde{y_i}=y_i(1-\phi_i+\omega\phi_i).
$$
It is easy to check that
$$
|1-\phi_i(t)+\omega\phi_i(t)|\leq 1 \qquad \big(t\in[0,1]\big),
$$
so, $\|\widetilde{y_i}\|\leq\|y_i\|$. This implies that
$\|\widetilde{y}\|=\|(\widetilde{y_1},\widetilde{y_2})\|\leq1$ and,
therefore,
\begin{equation}\label{example:norm(Id+omegaT)}
\|\Id+\omega T\|^2\geq\|\widetilde{y}+\omega
T(\widetilde{y})\|^2=\|\widetilde{y_1}+\omega
x^*(\widetilde{y})x_1\|^2+\|\widetilde{y_2}+\omega
x^*(\widetilde{y})x_2\|^2.
\end{equation}
Since, clearly,
$$
y_i-\widetilde{y_i}=(1-\omega)y_i\phi_i,\qquad (i=1,2)
$$
we deduce from \eqref{open-interval-property2},
\eqref{open-interval-property2-bis}, and the fact that
$\text{supp}(\phi_i)\subset \Delta_i$, that
$$
|x^*(y-\widetilde{y})|\leq
|x_1^*(y_1-\widetilde{y_1})|+|x_2^*(y_2-\widetilde{y_2})|\leq 4\,
\eps M,
$$
where $M=\max\{\|x_1^*\|,\|x_2^*\|\}$. Using this,
\eqref{open-interval-property}, and the fact that $\|x\|=1$, we
obtain the following for $i=1,2$:
\begin{align*}
\|\widetilde{y_i}+\omega\,
x^*(\widetilde{y})x_i\|&\geq\|\widetilde{y_i}+\omega\,
x^*(y)\,x_i\|-|x^*(y-\widetilde{y})|\\
 &\geq |\widetilde{y_i}(t_i)+\omega\,
x^*(y)\,x_i(t_i)|-4\,\eps M\\ & = |\omega\,y_i(t_i) +
\omega\,x^*(y)\,x_i(t_i)|-4\,\eps M\\
&= |y_i(t_i)+ x^*(y)\,x_i(t_i)|-4\,\eps M\\ &\geq
\|y_i+x^*(y)\,x_i\|-\eps-4\,\eps M.
\end{align*}
Finally, using this together with \eqref{example:norm(Id+T)} and
\eqref{example:norm(Id+omegaT)}, it is not hard to find a suitable
constant $K>0$ such that
$$
\|\Id+\omega T\|^2\geq\|\Id+T\|^2-K\,\eps,
$$
which finishes the proof.
\end{proof}

\vspace{0.5cm}

\noindent \textbf{\large $\bullet$} \textsc{\large Real case:}

The situation in the real case is far away from being so clear. On
the one hand, the proof of Lemma~\ref{lemma:4.3} remains valid if
the function $g$ is surjective (this substitutes Picard Theorem) and
then, the proofs of Theorems \ref{th:4.4} and \ref{th:4.5} are
valid. In addition, Example~\ref{example-new-property} was also
stated for the real case. The following result summarizes all these
facts.

\begin{theorem}\label{theorem-real-case}
Let $X$ be a real Banach space with dimension greater or equal than
two. Suppose that there exists a \textbf{surjective} analytic
function $g:\R\longrightarrow\R$ and a continuous function
$f:[|g(0),+\infty[\longrightarrow\R^+_0$, such that the norm
equality
$$
\|\Id+g(T)\|=f(\|g(T)\|)
$$
holds for every rank-one operator $T\in L(X)$.
\begin{enumerate}
\item[(a)] If $g(0)\neq -1/2$, then $X$ has the Daugavet property.
\item[(b)] If $g(0)=-1/2$, then the norm equality
$$
\|\Id - T\|=\|\Id + T\|
$$
holds for every rank-one operator $T\in L(X)$.
\item[(c)] The real space $X=C[0,1]\oplus_2 C[0,1]$ does not have the
Daugavet property but the norm equality
$$
\|\Id - T\|=\|\Id +T\|
$$
holds for every rank-one operator $T\in L(X)$.
\end{enumerate}
\end{theorem}

On the other hand, we do not know if a result similar to the above
theorem is true when the function $g$ is not onto. Let us give some
remarks about two easy cases:
$$
g(t)=t^2\quad (t\in\R) \qquad \text{and} \qquad g(t)=-t^2\quad
(t\in\R).
$$
In the first case, it is easy to see that if the norm equality
$$
\|\Id + T^2\|=f(\|T^2\|)
$$
holds for every rank-one operator, then $f(t)=1 + t$ and, therefore,
the interesting norm equality in this case is
\begin{equation}\label{eq:T2}
\|\Id +T^2\|=1+\|T^2\|.
\end{equation}
This equation is satisfied by every rank-one operator $T$ on a
Banach space $X$ with the Daugavet property. Let us also recall that
the equality
$$
|1+t^2|=1+|t^2|
$$
holds for every $t\in L(\R)\equiv\R$ (Remark~\ref{remark-dim-one}).
For the norm equality
$$
\|\Id -T^2\|=f(\|T^2\|),
$$
we are not able to get any information about the shape of the
function $f$. Going to the $1$-dimensional case, we get that
$$
|1-t^2|=\max\{1-|t^2|,|t^2|-1\}
$$
for every $t\in L(\R)\equiv\R$, but it is not possible that the
corresponding norm equality holds for all rank-one operators on a
Banach space with dimension greater than one (in this case,
$\|\Id-T^2\|\geq 1$). On the other hand, if a Banach space $X$ has
the Daugavet property, then
\begin{equation}\label{eq:-T2}
\|\Id-T^2\|=1+\|T^2\|
\end{equation}
for every rank-one operator $T\in L(X)$. Therefore, an interesting
norm equality of this form could be the above one.

Let us characterize the properties which flow out from the norm
equalities \eqref{eq:T2} and \eqref{eq:-T2}.

\begin{prop}
Let $X$ be a real Banach space.
\begin{enumerate}
\item[(a)] The following are equivalent:
\begin{equi}
\item $\|\Id+T^2\|=1+\|T^2\|$ for every rank-one operator $T$.
\item $\|\Id+x^*\otimes x\|=1+\|x^*\otimes x\|$ for $x^*\in
X^*$, $x\in X$ with $x^*(x)\geq0$.
\item For every $x\in S_X$, $x^*\in S_{X^*}$ with $x^*(x)\geq 0$,
and every $\eps>0$, there exists $y\in S_X$ such that
$$
\|x + y\| > 2-\eps \qquad \text{and} \qquad x^*(y)> 1-\eps.
$$
\end{equi}
\item[(b)] The following are equivalent:
\begin{equi}
\item $\|\Id - T^2\|=1+\|T^2\|$ for every rank-one operator $T$.
\item $\|\Id+x^*\otimes x\|=1+\|x^*\otimes x\|$ for $x^*\in
X^*$, $x\in X$ with $x^*(x)\leq 0$.
\item For every $x\in S_X$, $x^*\in S_{X^*}$ with $x^*(x)\leq 0$,
and every $\eps>0$, there exists $y\in S_X$ such that
$$
\|x + y\| > 2-\eps \qquad \text{and} \qquad x^*(y)> 1-\eps.
$$
\end{equi}
\end{enumerate}
\end{prop}

\begin{proof}
We give our arguments just for item (a), following (b) in the same
way.

$(i)\Rightarrow (ii)$. We consider $x\in X$ and $x^*\in X^*$ with
$x^*(x)>0$, and we write $T=x^*\otimes x$. Let us observe that the
rank-one operator $S=\big(x^*(x)\big)^{-1/2}\,T$ satisfies
$S^{\,2}=T$, so we get
$$
1+\|T\|=1+\|S^{\,2}\|=\|\Id+S^{\,2}\|=\|\Id+T\|.
$$
To finish the argument, it suffices to observe that any $x^*\otimes
x$ with $x^*(x)=0$ can be approximated in norm by operators of the
form $y^*\otimes y$ with $y^*(y)>0$, and that the set of rank-one
operators satisfying \eqref{DE} is closed.

$(ii)\Rightarrow(i)$. Just observe that for every rank-one operator
$T=x^*\otimes x$, it is clear that $T^2=x^*(x)\,T$ and, therefore,
$T^2=y^*\otimes x$ where $y^*=x^*(x)\,x^*$ with
$y^*(x)=\big(x^*(x)\big)^2\geq 0$.

Finally, for the equivalence between $(ii)$ and $(iii)$ just follow
the proof of \cite[Lemma~2.1]{KSSW} or \cite[Lemma~2.2]{WerSur}.
\end{proof}

\section{The new properties}\label{sec:newproperties}

Our aim in this section is to give some remarks concerning the
properties appearing at Theorems \ref{th:4.5} and
\ref{theorem-real-case}.b, i.e.\ we consider a non-trivial
multiplicative subgroup $A$ of $\T$ and study those Banach spaces
for which all rank-one operators satisfy the norm equality
$$
\|\Id + \omega\,T\|=\|\Id+T\|.
$$
In the real case, only one property arises; in the complex case,
there are infinitely many properties and we do not know if all of
them are equivalent.

Our first (easy) observation is that all these properties pass from
the dual of a Banach space to the space. We will see later that the
converse result is not valid.

\begin{remark}\label{remark:dualtospace}
{\slshape Let $X$ be a Banach space and let $\omega\in \T$. Suppose
that the norm equality
$$
\|\Id + \omega\,T\|=\|\Id+T\|
$$
holds for every rank-one operator $T\in L(X^*)$. Then, the same is
true for every rank-one operator on $X$.\ } Indeed, the result
follows routinely by just considering the adjoint operators of the
rank-one operators on $X$.
\end{remark}

Our next results deal with the shape of the unit ball of the Banach
spaces having any of the properties

By a \emph{slice} of a subset $A$ of a normed space $X$ we mean a
set of the form
$$
S(A,x^*,\alpha)=\bigl\{x\in A\ : \ \re x^*(x)>\sup \re x^*(A)-\alpha
\bigr\}
$$
where $x^*\in X^*$ and $\alpha\in \R^+$. If $X$ is a dual space, by
a \emph{weak$^*$-slice} of a subset $A$ of $X$ we mean a slice of
$A$ defined by a weak$^*$-continuous functional or, equivalently, a
weak$^*$-open slice $A$.

\begin{prop}
Let $X$ be a real or complex Banach space and let $A$ be a
non-trivial closed subgroup of $\T$. Suppose that the norm equality
$$
\|\Id + \omega\,T\|=\|\Id+T\|
$$
holds for every rank-one operator $T\in L(X)$ and every $\omega\in
A$. Then, the slices of $B_X$ and the weak$^*$-slices of $B_{X^*}$
have diameter greater or equal than
$$
2-\inf \bigl\{|1+\omega|\ : \ \omega\in A\bigr\}.
$$
\end{prop}

\begin{proof} We give the arguments only for slices of $B_X$, being
the proof for weak$^*$-slices of $B_{X^*}$ completely analogous. We
fix $x^*\in S_{X^*}$ and $0<\alpha<2$. Given $0<\eps<\alpha$, we
take $x\in S_{X}$ such that $\re x^*(x)>1-\eps$ and so, in
particular, $x\in S(B_X,x^*,\alpha)$. We define the rank-one
operator $T=x^*\otimes x$ and observe that
$$
\|\Id+T\|\geq\|x + x^*(x)\,x\|=\|x\|\,|1+x^*(x)|\geq|1+\re
x^*(x)|>2-\eps.
$$
By hypothesis, for every $\omega\in A$ we may find $y\in S_X$ such
that
$$
\|y+\omega\,x^*(y)\,x\|>2-\eps
$$
so, in particular,
$$
|x^*(y)|>1-\eps.
$$
We take $\xi\in \T$ such that $\xi\,x^*(y)=|x^*(y)|$ and we deduce
that $$\xi\,y\in S(B_X,x^*,\alpha).$$ From the inequalities
\begin{align*}
\|\xi y-x\| &=\big\|\xi y +\omega x^*(\xi y) x-\omega x^*(\xi y) x -x\big\|\\
&\geq \big\|\xi y +\omega x^*(\xi y) x\big\| -\big|1+\omega x^*(\xi y)\big| \\
&> 2-\eps -\big|1+\omega |x^*(y)|\big|
\end{align*}
and
\begin{align*}
\big|1+\omega |x^*(y)|\big|&\leq\big|1+\omega
+\omega(|x^*(y)|-1)\big|\\ &\leq |1+\omega|+
\big|1-|x^*(y)|\big|<|1+\omega|+\eps,
\end{align*}
we get
\begin{equation*}
\|\xi y-x\|\geq 2-2\,\eps - |1+\omega|.\qedhere
\end{equation*}
\end{proof}

It is well-known that the unit ball of a Banach space $X$ with the
Radon-Nikod\'{y}m property has many denting points and the unit ball of
the dual of an Asplund space has many weak$^*$-denting points.
Recall that $x_0\in B_X$ is said to be a \emph{denting point} of
$B_X$ if it belongs to slices of $B_X$ with arbitrarily small
diameter. If $X$ is a dual space and the slices can be taken to be
weak$^*$-open, then we say that $x_0$ is a \emph{weak$^*$-denting
point}. We refer to \cite{Bou,D-U} for more information on these
concepts.

\begin{corollary}\label{cor-newpropertynodenting}
Let $X$ be a real or complex Banach space and let $\omega\in
\T\setminus\{1\}$. If the norm equality
$$
\|\Id + \omega\,T\|=\|\Id+T\|
$$
holds for every rank-one operator $T\in L(X)$, then $B_X$ does not
have any denting point and $B_{X^*}$ does not have any $w^*$-denting
point. In particular, $X$ is not an Asplund space and it does not
have the Radon-Nikod\'{y}m property.
\end{corollary}

In view of this result, it is easy to show that the converse of
Remark~\ref{remark:dualtospace} is not true.

\begin{example}
{\slshape The real or complex space $X=C[0,1]$ satisfies that the
norm equality
$$
\|\Id+\omega\,T\|=\|\Id+T\|
$$
holds for every rank-one operator $T\in L(X)$ and every $\omega\in
\T$, in spite of the fact that for every $\omega\in \T\setminus
\{1\}$, there is a rank-one operator $S\in L(X^*)$ such that}
$$
\|\Id + \omega\, S\|<\|\Id + S\|.
$$
Indeed, the first assertion follows from the fact that $X$ has the
Daugavet property; the second one follows from
Corollary~\ref{cor-newpropertynodenting} since the unit ball of
$X^*$ is plenty of denting points.
\end{example}

One of the properties we are dealing with is related to a property
for one-codimensional projections.

\begin{remark}
Consider a Banach space $X$ such that the equality
$$\|\Id + T\|=
\|\Id - T\|
$$
holds true for every rank-one operator $T\in L(X)$. If we apply this
equality to an operator $P$ which is a rank-one projection, we get
$$
\|\Id - P\| \geqslant 2,
$$
i.e. every one-codimensional projection in $L(X)$ is at least of
norm $2$. Such spaces were introduced recently \cite{IvKa} and are
called ``spaces with bad projections''.
\end{remark}

\end{document}